\documentclass{amsart}

\usepackage[T1]{fontenc}
\usepackage[latin1]{inputenc}
\usepackage{amscd}
\usepackage{epsfig}
\makeatletter

\theoremstyle{plain}    
\newtheorem{thm}{Theorem}[section]
\newtheorem{defn}[thm]{Definition}
\newtheorem{question}[thm]{Question}
\numberwithin{equation}{section} %% Comment out for sequentially-numbered
\numberwithin{figure}{section} %% Comment out for sequentially-numbered
\theoremstyle{plain}    
\newtheorem{cor}[thm]{Corollary} %%Delete [thm] to re-start numbering
\newtheorem{lem}[thm]{Lemma} %%Delete [thm] to re-start numbering
\theoremstyle{plain}    
\newtheorem{prop}[thm]{Proposition} %%Delete [thm] to re-start numbering
\theoremstyle{remark}
\newtheorem{rem}{Remark}[thm]
\newtheorem{warning}[thm]{Warning}
\newtheorem{construction}[thm]{Construction}
\theoremstyle{remark}    
\newtheorem{notation}[thm]{Notation}

\def\ringB{R}
\newcommand{\J}{\mbox{$\mathcal{J}$}} 
\newcommand{\inter}{{\ell_1.\ell_2}}
 
\def\factor#1.#2.{\left. \raise 2pt\hbox{$#1$} \right/\hskip -2pt\raise -2pt\hbox{$#2$}}
\def\ring#1.{\Cal O_{#1}} \newcommand\Cal{\mathcal}
\usepackage[arrow,matrix,curve,ps]{xy}
\newcommand{\C}{\mbox{$\mathbb{C}$}}

\renewcommand{\O}{\mbox{$\mathcal{O}$}}
\renewcommand{\P}{\mathbb{P}}

\newcommand{\ZZ}{\mbox{$\mathbb{Z}$}}

\newcommand{\sheafspec}{\mbox{\bf Spec}}

\DeclareMathOperator{\Aut}{Aut}

\DeclareMathOperator{\ptRC}{{RC_{\bullet}}}
\DeclareMathOperator{\twoptRC}{{RC^2_{\bullet}}}

\DeclareMathOperator{\Univrc}{{Univ^{rc}}}

\DeclareMathOperator{\ptUnivrc}{{Univ^{rc}_{\bullet}}}
\DeclareMathOperator{\twoptUnivrc}{{Univ_{\bullet}^{rc, 2}}}
\DeclareMathOperator{\charact}{char}

\DeclareMathOperator{\codim}{codim}
\DeclareMathOperator{\Dubbies}{Dubbies^n}
\DeclareMathOperator{\Ext}{Ext}

\DeclareMathOperator{\Hom}{Hom}
\DeclareMathOperator{\Hombir}{Hom_{bir}^n}
\DeclareMathOperator{\Image}{Image}

\DeclareMathOperator{\Num}{Num}
\DeclareMathOperator{\Pic}{Pic}

\DeclareMathOperator{\RatCurves}{RatCurves^n}
\DeclareMathOperator{\RC}{RatCurves^n}
\DeclareMathOperator{\red}{red}

\DeclareMathOperator{\Spec}{Spec}

\DeclareMathOperator{\im}{Im}

\DeclareMathOperator{\psl}{PGL} 
\def\twopttau{\tau^{\operatorname{rc}, 2}}

\makeatother

\begin{document}

\title[Are rational curves determined by tangent vectors?]{Are minimal
  degree rational curves determined by their tangent vectors?}

\begin{abstract}
  Let $X$ be a projective variety which is covered by rational curves,
  for instance a Fano manifold over the complex numbers. In this
  setup, characterization and classification problems lead to the
  natural question: ``Given two points on $X$, how many minimal degree
  rational curve are there which contain those points?''. A recent
  answer to this question led to a number of new results in
  classification theory. As an infinitesimal analogue, we ask ``How
  many minimal degree rational curves exist which contain a prescribed
  tangent vector?''
  
  In this paper, we give sufficient conditions which guarantee that
  every tangent vector at a general point of $X$ is contained in at
  most one rational curve minimal degree. As an immediate application,
  we obtain irreducibility criteria for the space of minimal rational
  curves.
\end{abstract}

\date{\today}

\author{Stefan Kebekus and S\'andor J.~Kov\'acs}



\thanks{S.~Kebekus gratefully acknowledges support by the
  Forschungsschwerpunkt ``Globale Methoden in der komplexen
  Geometrie'' of the Deutsche Forschungsgemeinschaft. S.~Kov\'acs was
  supported in part by NSF Grants DMS-0196072 and DMS-0092165 and by a
  Sloan Research Fellowship.}

\address{Stefan Kebekus, Mathematik VIII, Universit\"at Bayreuth, 95440
  Bayreuth, Germany}
\email{{stefan.kebekus@uni-bayreuth.de}}
\urladdr{{http://btm8x5.mat.uni-bayreuth.de/$\sim$kebekus}}

\address{\noindent
S\'andor Kov\'acs, University of Washington, Department of
  Mathematics, Box 354350, Seattle, WA 98195, U.S.A.}  
\email{{kovacs@math.washington.edu}}
\urladdr{{http://www.math.washington.edu/$\sim$kovacs}}

\maketitle
\tableofcontents

\section{Introduction}

The study of rational curves of minimal degree has proven to be a very
useful tool in Fano geometry. The spectrum of application covers
diverse topics such as deformation rigidity, stability of the tangent
sheaf, classification problems or the existence of non-trivial finite
morphisms between Fano manifolds; see \cite{Hwa} for an overview.

In this paper we will consider the situation where $X$ is a projective
variety, which is covered by rational curves, e.g.~a Fano manifold
over $\C$. An example of that is $\P^n$, which is covered by lines.
The key point of many applications of minimal degree rational curves
is showing that the curves in question are similar to lines in certain
respects. For instance, one may ask:

\begin{question} \label{q:iota}
  Under what conditions does there exist a unique minimal degree
  rational curve containing two given points?
\end{question}

This question found a sharp answer in \cite{Keb00a}, see \cite{CMS00}
and \cite{Keb01b} for a number of applications.  The argument used
there is based on a criterion of Miyaoka, who was the first to observe
that if the answer to the question is ``No'', then a lot of minimal
degree curves are singular. We refer to \cite[Prop.~V.3.7.5]{K96} for
a precise statement.

As an infinitesimal analogue of this question one may ask the following:

\begin{question} \label{q:tau}
  Are there natural conditions that guarantee that a minimal degree
  rational curve is uniquely determined by a tangent vector?
\end{question}

Although a definite answer to the latter question would be as
interesting as one to the former, it seems that Question~\ref{q:tau}
has hardly been studied before.  This paper is a first attempt in that
direction. We give a criterion which parallels Miyaoka's approach.

\begin{thm}\label{thm:main-criterion1}
  Let $X$ be a projective variety over an algebraically closed field
  $k$ and $H \subset \RatCurves(X)$ a proper, covering family of
  rational curves such that none of the associated curves has a
  cuspidal singularity. If $\charact(k) \not = 0$, assume additionally
  that there exists an ample line bundle $L \in \Pic(X)$ such that for
  every $\ell \in H$ the intersection number $L.\ell$ of $L$ is
  coprime to $\charact(k)$.
  
  Then, if $x \in X$ is a general point, all curves associated with
  the closed subfamily
  $$
  H_x := \{ \ell \in H \,|\, x \in \ell\} \subset H
  $$
  are smooth at $x$ and no two of them share a common tangent
  direction at $x$.
\end{thm}

\begin{rem}
  In Theorem~\ref{thm:main-criterion1} we do not assume that $H$ is
  irreducible or connected. That will later be important for the
  applications.
\end{rem}

\begin{rem}
  We refer the reader to Chapter~\ref{sec:ratcurves-review} for a
  brief review of the space $\RatCurves(X)$ of rational curves. The
  volume \cite{K96} contains a thorough discussion.
  
  If $H \subset \RatCurves(X)$ is an irreducible component, it is
  known that $H$ is proper if there exists a line bundle $L \in
  \Pic(X)$ that intersects a curve $\ell \in H$ with multiplicity
  $L.\ell =1$.
\end{rem}

For complex projective manifolds we give another result. To formulate
the setup properly, pick an irreducible component $H \subset
\RatCurves(X)$ such that
\begin{enumerate}
\item the rational curves associated with $H$ dominate $X$,
  
\item for a general point $x\in $, the closed subfamily $H_x$ is
  proper.
\end{enumerate}
Let $\tilde U$ be the universal family, which is a $\P^1$-bundle over
$H$. The tangent map of the natural projection $\iota: \tilde U\to X$,
restricted to the relative tangent sheaf $T_{\tilde U/H}$, gives rise
to a rational map $\tau$:
$$ 
\xymatrix{ & & {\P(T_X^\vee)} \ar[d] \\ 
{\tilde U} \ar@{-->}@/^0.5cm/[rru]^{\tau} 
\ar[rr]^{\iota}_{\txt{\scriptsize evaluation}} 
\ar[d]_{\pi}^{\txt{\scriptsize $\P^1$-bundle}} & & {X} \\
{H} }
$$
It has been shown in \cite{Keb00a} that $\tau$ is generically
finite. Examples of rationally connected manifolds, however, seem to
suggest that the tangent map $\tau$ is generically injective for a
large class of varieties. Our main result supports this claim.

\begin{thm}\label{thm:main-criterion2}
  Let $X$ be a smooth projective variety over the field of complex
  numbers and let $H \subset \RatCurves(X)$ be the union of
  irreducible components such that the subfamily $H_x$ is proper for
  all points $x\in X$, outside a subvariety $S \subset X$ of
  codimension at least~$2$. 

  Then $\tau$ is generically injective, unless the curves associated
  with the closed subfamily $H^{\text{cusp}} \subset H$ of cuspidal
  curves dominate $X$, and the subvariety
  $$
  D := \{ x\in X\, |\, \text{ $\exists \ell \in H^{\text{cusp}}:
    \ell$ has a cuspidal singularity at $x$} \},
  $$
  where curves have cuspidal singularities, has codimension 1.
\end{thm}

\begin{rem}
  It is known that the family $H_x$ is proper for a general point $x
  \in X$ if $H$ is a ``maximal dominating family of rational curves of
  minimal degrees'', i.e., if the degrees of the curves associated
  with $H$ are minimal among all irreducible components of
  $\RatCurves(X)$ which satisfy condition~(1) from above.
  
  The assumption that $H_x$ is proper for all points outside a set of
  codimension 2, however, is restrictive.
\end{rem}

Although we consider these results only a first step toward a
satisfactory answer to Question~\ref{q:tau}, we are already able to
present several applications in Section~\ref{sect:applications}.

\subsubsection*{Acknowledgements}

Parts of this paper have been worked out while the first named author
visited the University of Washington at Seattle, the University of
British Columbia at Vancouver and Princeton University as well as
while the second named author visited the Isaac Newton Institute for
Mathematical Sciences at Cambridge. Both authors are grateful to these
institutions for their hospitality. S.~Kebekus would like to thank
K.~Behrend and J.~Koll\'ar for the invitations and for numerous
discussions. S.~Kov\'acs would like to thank the organizers, Alessio
Corti, Mark Gross, and Miles Reid for the invitation to the 'Higher
Dimensional Geometry Programme' of the Newton Institute.

After the main part of this paper was written, J.-M.~Hwang has
informed us
%the authors 
that, together with N.~Mok, they have shown a
statement similar to, but somewhat stronger than
Theorem~\ref{thm:main-criterion2}. Their unpublished proof uses
entirely different methods. To the best of our knowledge, there is no
other result similar to Theorem~\ref{thm:main-criterion1}.

\section{$\P^1$-bundles with double sections}
\label{sec:P1_bundles}

This preliminary section discusses $\P^1$-bundles with an irreducible
double section. Most results here are fairly elementary.  We have,
however, chosen to include detailed proofs for lack of a suitable
reference.

Throughout the present section let $\lambda : \Lambda \to B$ be a
$\P^1$-bundle over a normal variety $B$, i.e., a morphism whose
scheme-theoretic fibers are all isomorphic to $\P^1$. Note that we do
not assume here that $\lambda$ is trivial in the Zariski topology. Let
$\sigma:B\to \Lambda$ be a section of $\lambda$, $\Sigma_{\red}
=\sigma(B)_{\red} \subset \Lambda$, and let $\Sigma \subset \Lambda$
be the first infinitesimal neighborhood of $\Sigma_{\red}$ in
$\Lambda$. That is, if $\Sigma_{\red}$ is defined by the sheaf of
ideals $\J=\O_\Lambda(-\Sigma_{\red})$, then $\Sigma$ is defined by
the sheaf $\J^2$.  Our aim is to relate properties of $\Lambda$ with
those of its subscheme $\Sigma$.

\subsection{The Picard group of the double section}
Recall from \cite[III.~Ex.4.6]{Ha77} qthat there exists a short exact
sequence of sheaves of Abelian groups, sometimes called that
``truncated exponential sequence'' in the literature (eg.
\cite[sect.~2]{BBI00})
\begin{equation}
  \label{eq:truncated_exp} 
  0 \longrightarrow \underbrace{\J \big/\J^2}_{=
    N_{\Sigma_{\red}|\Lambda}^{\vee}}
  \overset{\alpha}{\longrightarrow} \O_{\empty\Sigma}^*
  \overset{\beta}{\longrightarrow} \O_{\empty\Sigma_{\red}}^*
  \longrightarrow 1.
\end{equation}
Here $N_{\Sigma_{\red}|\Lambda}^{\vee}$ is the conormal bundle,
$\beta$ is the canonical restriction map and $\alpha$ is given by
$$
\begin{array}{cccc}
  \alpha: & (\J \big/ \J^2, +) &\to & \O_{\empty\Sigma}^* \\ 
  & f &\mapsto & 1+f. 
\end{array}
$$ In our setup, where $\Sigma_{\red} \simeq B$ is a section, the
truncated exponential sequence~\eqref{eq:truncated_exp} is canonically
split. Locally we can write the splitting as follows.  Assume that we
are given an affine open subset $U_\alpha \subset \Sigma$ and an
invertible function $f_\alpha \in
\O_{\empty\Sigma}^*(U_\alpha)$. Then, after shrinking $U_\alpha$, if
needed, we will find a {bundle} coordinate $y_a$, identify
$$
\O_{\empty\Sigma}^*(U_\alpha) \simeq
[ \O_{\empty\Sigma_{\red}}(U_\alpha) \otimes k[y_\alpha]/(y_\alpha^2)]^*
$$
and write accordingly
$$
f_\alpha = g_\alpha+h_\alpha\cdot y_\alpha
$$
where $g_\alpha \in \O^*_{\Sigma_{\red}}(U_\alpha)$ and $h_\alpha
\in \O_{\Sigma_{\red}}(U_\alpha)$. With this notation, the splitting
of sequence~\eqref{eq:truncated_exp} decomposes $f_\alpha$ as
$$ f_\alpha = g_\alpha \cdot\underbrace{\left[
1+\frac{h_\alpha}{g_\alpha}\cdot y_\alpha \right]}_{\in
\im\left(\alpha_{U_\alpha}\right)}
$$
As a direct corollary to the splitting of~\eqref{eq:truncated_exp}
we obtain a canonical decomposition of the Picard group
\begin{equation}
  \label{eq:splitting_of_Pic}
  \Pic(\Sigma) = \Pic(\Sigma_{\red}) \times
  H^1(\Sigma_{\red},N_{\Sigma_{\red}|\Lambda}^{\vee}).
\end{equation}

\subsection{The cohomology class of a line bundle}

Let $L \in \Pic(\Lambda)$ be a line bundle. Using the
decomposition~\eqref{eq:splitting_of_Pic} from above, we can associate
to $L$ a class $c(L) \in H^1(\Sigma_{\red},
N_{\Sigma_{\red}|\Lambda}^{\vee})$. As this class will be important
soon, we will now find a {\v C}ech-cocycle in $Z^1(U_\alpha,
N_{\Sigma_{\red}|\Lambda}^{\vee})$ that represents $c(L)$.

To this end, find a suitable open affine cover $U_\alpha$ of $\Sigma$
such that $L|_{U_\alpha}$ is trivial for all $\alpha$ and where
{bundle} coordinates $y_\alpha$ exist.  Let $f_\alpha \in L(U_\alpha)$
be a collection of nowhere vanishing sections which we write in local
coordinates as $f_\alpha = g_\alpha + h_\alpha\cdot y_\alpha$. Using
the $U_\alpha$-coordinates on the intersection $U_\alpha \cap
U_\beta$, the transition functions for the line bundle are thus
written as
$$ \frac{f_\alpha}{f_\beta} = \frac{g_\alpha + h_\alpha\cdot
y_\alpha}{g_\beta + h_\beta\cdot y_\alpha} = \frac{g_\alpha}{g_\beta}
\cdot \left[ 1 + \left( \frac{h_\alpha}{g_\alpha} -
\frac{h_\beta}{g_\beta} \right)y_\alpha \right] \in
\O^*_{\Sigma}(U_{\alpha\beta})
$$ 
In other words, the class of $c(L) \in H^1(\Sigma_{\red},
N_{\Sigma_{\red}|\Lambda}^{\vee})$ is represented by the {\v C}ech
cocycle
\begin{equation}
  \label{eq:cocycle}
  \left( \frac{h_\alpha}{g_\alpha} - \frac{h_\beta}{g_\beta}
  \right)y_\alpha \in Z^1(U_\alpha, N_{\Sigma_{\red}|\Lambda}^{\vee})  
\end{equation}

\subsection{Vector bundle sequences associated to line bundles}

Consider the ideal sheaf sequence for $\Sigma_{\red} \subset \Sigma$.
$$
0 \longrightarrow \J \big/\J^2 \longrightarrow \O_{\empty\Sigma}
\longrightarrow \O_{\empty\Sigma_{\red}} \longrightarrow 0
$$

\begin{warning}
  It should be noted that $\Sigma_{\red}$ is not a Cartier-divisor in
  $\Sigma$ since its ideal sheaf, $\J \big/\J^2 \simeq
  N_{\Sigma_{\red}|\Lambda}^{\vee}$ is {\em not} a locally free
  $\O_{\Sigma}$-module. Furthermore, the restriction of the ideal
  sheaf of $\Sigma_{\red}$ in $\Lambda$ to $\Sigma$, $\J\otimes
  \O_{\Sigma} \simeq\J\big/\J^3$, is not isomorphic to the ideal sheaf
  of $\Sigma_{\red}$ in $\Sigma$, $\J \big/\J^2 \not \simeq
%%  \O_{\Lambda}(-\Sigma_{\red})
\J\otimes \O_{\Sigma}.$ In fact, $\J\otimes \O_{\Sigma}$ is not even a
subsheaf of $\O_{\Sigma}$.
\end{warning}

\begin{construction}\label{const:VB}
 Let $L \in \Pic(\Sigma)$ be a line bundle.  By abuse of
notation, identify $\Sigma_{\red}$ with $B$ and consider
$L|_{\Sigma_{\red}}$ a line bundle on $B$. Then twist the above
sequence with the locally free $\O_{\Sigma}$-module $L \otimes
\lambda^*(L^\vee|_{\Sigma_{\red}})$, and obtain the following sequence
of $\O_{\Sigma}$-modules,
\begin{equation}
  \label{eq:seq1e}
  0 \longrightarrow N_{\Sigma_{\red}|\Lambda}^{\vee} \longrightarrow
  L\otimes \lambda^*(L^\vee|_{\Sigma_{\red}}) \longrightarrow
  \O_{\Sigma_{\red}} \longrightarrow 0.
\end{equation}
Finally, consider the push-forward to $B$:
\begin{equation}
  \label{eq:push_second_seq}
  0 \longrightarrow N_{\Sigma_{\red}|\Lambda}^{\vee} \longrightarrow 
  \underbrace{\lambda_*(L)\otimes L^\vee|_{\Sigma_{\red}}}_{=:
    \mathcal E_L} \overset{A}{\longrightarrow} \O_{\Sigma_{\red}}
  \longrightarrow 0.
\end{equation}

We obtain a vector bundle $\mathcal E_L$ of rank two on $B$ which is
presented as an extension of two line bundles. The surjective map
$\mathcal E_L\to \O_B$ induces a section $\sigma_L:B\to \P(\mathcal
E)$. We will use this notation later and also extend it to line
bundles, $L\in\Pic(\Lambda)$, by $\mathcal E_L:=\mathcal
E_{L|_{\Sigma}}$ and $\sigma_L=\sigma_{L|_{\Sigma}}$.  Observe that
$(\P(\mathcal E_L), \sigma_L)$ depends on $L$ only up to a twist by a
line bundle pulled back from $B$. I.e., for $M\in \Pic(B)$, $\mathcal
E_{L\otimes \lambda^*M}\simeq\mathcal E_L$ and $\sigma_{L\otimes
\lambda^*M}=\sigma_L$.
\end{construction}

Much of our further argumentation is based on the following
observation.

\begin{prop}
  \label{prop:PicAndExt}
  Let $L \in \Pic(\Sigma)$ be a line bundle and $c(L) \in
  H^1(\Sigma_{\red}, N_{\Sigma_{\red}|\Lambda}^{\vee})$ the class
  defined above. Then $c(L)$ coincides with the extension class
  $$
  e(L) \in \Ext^1(\O_{\Sigma_{\red}},
  N_{\Sigma_{\red}|\Lambda}^{\vee}) = H^1(\Sigma_{\red},
  N_{\Sigma_{\red}|\Lambda}^{\vee} )
  $$
  of the vector bundle sequence~\eqref{eq:push_second_seq}. 
In particular, the map
  $$
  \begin{array}{rccc}
    e : & (\Pic(\Sigma),\otimes) & \to & (H^1(\Sigma_{\red},
    N_{\Sigma_{\red}|\Lambda}^{\vee} ),+) \\ & L & \mapsto &
    \text{extension class of sequence~\eqref{eq:push_second_seq}}
  \end{array}
  $$
  is a homomorphism of groups.
\end{prop}

\begin{proof}
  The proof relies on an explicit calculation in {\v C}ech cohomology.
  We will choose a sufficiently fine cover $U_\alpha$ of
  $\Sigma_{\red}$ and produce a {\v C}ech cocycle in $Z^1(U_\alpha,
  N_{\Sigma_{\red}|\Lambda}^{\vee})$ that represents the extension
  class $e(L)$. It will turn out that this cocycle equals the one that
  we have calculated in~\eqref{eq:cocycle} above for $c(L)$.
  
  We keep the notation from above and let $f_\alpha \in L(U_\alpha)$
  be a collection of nowhere-vanishing sections of $L$.  Such sections
  can be naturally seen to give local splittings of the
  sequences~\eqref{eq:seq1e} and~\eqref{eq:push_second_seq}.
  Explicitly, if we write $f_\alpha = g_\alpha + h_\alpha \cdot
  y_\alpha$, then
  $$
  \frac{f_\alpha}{g_\alpha} = 1+\frac{h_\alpha}{g_\alpha}\cdot
  y_\alpha \in (L \otimes L^\vee|_{\Sigma_{\red}})(U_\alpha)
  $$
  are nowhere-vanishing sections of $L \otimes
  L^\vee|_{\Sigma_{\red}}$ and the splitting takes the form
  $$
  \begin{array}{rccc}
    s_\alpha : & \O_{\Sigma_{\red}}(U_\alpha) &\to& (L \otimes
    L^\vee|_{\Sigma_{\red}})(U_\alpha) \\
    & 1 & \mapsto & 1+\frac{h_\alpha}{g_\alpha}\cdot y_\alpha
  \end{array}
  $$
  By construction of $\Ext^1$, we obtain the extension class as the
  homology class represented by the {\v C}ech cocycle
  $$
  s_\alpha(1) - s_\beta(1) \in \ker(A)(U_{\alpha\beta}) \simeq
  N_{\Sigma_{\red}|\Lambda}^{\vee} (U_{\alpha\beta})
  $$
  This difference is given by the following section in
  $N_{\Sigma_{\red}|\Lambda}^{\vee} (U_{\alpha\beta})$ which yields
  the required cocycle.
  $$
  \left( 1+\frac{h_\alpha}{g_\alpha}\cdot y_\alpha \right) -\left(
    1+\frac{h_\beta}{g_\beta}\cdot y_\alpha \right) = \left(
    \frac{h_\alpha}{g_\alpha} - \frac{h_\beta}{g_\beta}
  \right)y_\alpha \in Z^1(U_\alpha, N_{\Sigma_{\red}|\Lambda}^{\vee})
  $$ That, however, is the same cocycle which we have obtained above
  in formula~\eqref{eq:cocycle} for the class $c(L)$. The proof of
  Proposition~\ref{prop:PicAndExt} is therefore finished.
\end{proof}

\subsection{The reconstruction of the $\P^1$-bundle from a double section}

It is a remarkable fact that the restriction of an ample line bundle
$L \in \Pic(\Lambda)$ to a double section carries enough information
so that the whole $\P^1$-bundle $\Lambda$ can be reconstructed. The
proof is little more than a straightforward application of
Proposition~\ref{prop:PicAndExt}. We are grateful to Ivo Radloff who
showed us how to use extension classes to simplify our original proof.

\begin{notation}
  Let $(\Lambda, \sigma)$ and $(\Lambda', \sigma')$ be two
  $\P^1$-bundles with sections over $B$. We say that $(\Lambda,
  \sigma)$ and $(\Lambda', \sigma')$ are \emph{isomorphic pairs (over
  $B$)} if there exists a morphism $\gamma : \Lambda \to \Lambda'$, an
  \emph{isomorphism of pairs}, such that $\gamma$ is a $B$-isomorphism
  of $\P^1$-bundles and $\gamma\circ\sigma = \sigma'$.  Sometimes we
  will refer to these pairs by the image of the section: $(\Lambda,
  \sigma(B))$, in which case the meaning of \emph{isomorphic pairs}
  should be clear.
\end{notation}

\begin{thm}\label{thm:isomorphism_of_pairs}
  Given a line bundle $L \in \Pic(\Lambda)$, which is not the
  pull-back of a line bundle on $B$, let $\mathcal E_L$ and $\sigma_L$
  be as in \ref{const:VB}.  Consider the relative degree $d \in
  \mathbb Z\setminus \{0\}$ of $L$, i.e., the intersection number with
  fibers of $\lambda$.  If $d$ is coprime to $\charact(k)$, then
  $(\Lambda, \sigma)$ and $(\P(\mathcal E_L), \sigma_L)$ are
  isomorphic pairs over $B$.
\end{thm}

\begin{proof}
  Let $H := \O_\Lambda(\Sigma_{\red})= \J^\vee$. Then
  $\Lambda\simeq \P(\lambda_*H)$ and $\sigma : B\to\Lambda$ is the
  section associated to the surjection, $\lambda_*H\to
  \lambda_*(H|_{\Sigma_{\red}})$.

First we would like to prove that
$\lambda_*H\simeq\lambda_*(H|_{\Sigma})$. Indeed, consider the
sequence,
$$
0\longrightarrow H\otimes \J^2 \simeq \J \longrightarrow 
H \longrightarrow H|_{\Sigma} \longrightarrow 0.
$$ We need to prove that $\lambda_* \J\simeq R^1\lambda_* \J\simeq
0$. However, that follows from considering the pushforward of the
sequence,
$$
0\longrightarrow  \J \longrightarrow 
\O_\Lambda \longrightarrow \O_{\Sigma_{\red}} \longrightarrow 0,
$$ since
$\lambda_*\O_\Lambda\simeq\lambda_*\O_{\Sigma_{\red}}\simeq\O_B$, and
$R^1\lambda_*\O_\Lambda\simeq 0$.

This implies the statement for $L=H$, that is, we obtain that
  $(\Lambda, \sigma)$ and $(\P(\mathcal E_H), \sigma_H)$ are
isomorphic pairs over $B$ (cf.~\cite[II.7.9]{Ha77}).

In order to finish the proof, we are going to prove that $(\P(\mathcal
  E_H), \sigma_H)$ and $(\P(\mathcal E_L), \sigma_L)$ are isomorphic
  pairs over $B$ for any $L\in \Pic(\Lambda)$.  In fact, it suffices
  to show that the extension classes of the following sequences are
  the same up to a non-zero scalar multiple.
\begin{equation}
\begin{aligned}
  \label{eq:HandL}
  0 \longrightarrow N_{\Sigma_{\red}|\Lambda}^{\vee} \longrightarrow
  \lambda_* (H|_\Sigma) &\otimes H^\vee|_{\Sigma_{\red}}
  \longrightarrow \O_{\Sigma_{\red}} \longrightarrow 0\\ 0
  \longrightarrow N_{\Sigma_{\red}|\Lambda}^{\vee} \longrightarrow
  \lambda_* (L|_\Sigma) &\otimes L^\vee|_{\Sigma_{\red}}
  \longrightarrow \O_{\Sigma_{\red}} \longrightarrow 0.
\end{aligned}
\end{equation}

  Recall that $\Pic(\Lambda) = {\mathbb Z} \times \Pic(B)$ so that we
  can write $L \in H^{\otimes d}\otimes \lambda^*M$ for an appropriate
  $M\in\Pic(B)$. By Proposition~\ref{prop:PicAndExt} this implies that
  the extension classes of the sequences~\eqref{eq:HandL} are given by
  $c\left(H|_\Sigma\right)$ and $c\left(H|_\Sigma^{\otimes d}\right) =
  d\cdot c\left(H|_\Sigma\right)$. In particular, they differ only by
  the non-zero factor $d \in k$.
\end{proof}

\begin{warning}
  The construction of the vector bundle $\mathcal E_L$ and
  Proposition~\ref{prop:PicAndExt} use only the restriction
  $L|_\Sigma$. It may thus appear that
  Theorem~\ref{thm:isomorphism_of_pairs} could be true without the
  assumption that $L \in \Pic(\Lambda)$ and that one could allow
  arbitrary line bundles $L \in \Pic(\Sigma)$ instead. That, however,
  is wrong and counterexamples do exist. Note that the proof of
  Theorem~\ref{thm:isomorphism_of_pairs} uses the fact that $L$ is
  contained in $\mathbb Z\times\Pic(B)$ which is not true in general
  if $L \in \Pic(\Sigma)$ is arbitrary.
\end{warning}

The assumption that $d$ be coprime to $\charact(k)$ is actually
necessary in Theorem~\ref{thm:isomorphism_of_pairs}, as shown by the
following simple corollary of Proposition~\ref{prop:PicAndExt} and of
the proof of Theorem~\ref{thm:isomorphism_of_pairs}.

\begin{cor}\label{cor:divisible}
  Using the same notation as in
  Theorem~\ref{thm:isomorphism_of_pairs}, assume that $d$ is divisible
  by $\charact(k)$. Then
  $$
  \lambda_* (L|_\Sigma) \otimes L^\vee|_{\Sigma_{\red}}
  \simeq N_{\Sigma_{\red}|\Lambda}^{\vee}\oplus \O_{\Sigma_{\red}}.
  $$
\end{cor}

%%% Local Variables: 
%%% mode: latex
%%% TeX-master: "Birationality_Thms"
%%% End: 

\section{Dubbies}

Throughout the proofs of Theorems~\ref{thm:main-criterion1} and
\ref{thm:main-criterion2}, which we give in
Sections~\ref{chap:proof11} and \ref{chap:proof12} below, we will
assume that $X$ contains pairs of minimal rational curves which
intersect tangentially in at least one point. A detailed study of
these pairs and their parameter spaces will be given in the present
chapter. The simplest configuration is the following:

\begin{defn}
  A \emph{dubby} is a reduced, reducible curve, isomorphic to the
  union of a line and a smooth conic in $\P^2$ intersecting
  tangentially in a single point.
\end{defn}

\begin{rem}
  The definition may suggest at first glance that one component of a
  dubby is special in that it has a higher degree than the other. We
  remark that this is not so. A dubby does not come with a natural
  polarization. In fact, there exists an involution in the automorphism
  group that swaps the irreducible components.
\end{rem}

Later we will need the following estimate for the dimension of the
space of global sections of a line bundle on a dubby. Let $\ell =
\ell_1 \cup \ell_2$ be a dubby and $L \in \Pic(\ell)$ a line bundle.
We say that $L$ has type $(d_1,d_2)$ if the restrictions of $L$ to the
irreducible components $\ell_1$ and $\ell_2$ have degree $d_1$ and
$d_2$, respectively.

\begin{lem}\label{lem:h0_of_dubbies}
  Let $\ell$ be a dubby and $L \in \Pic(\ell)$ a line bundle of type
  $(d_1, d_2)$. Then $h^0(\ell, L)\geq d_1+d_2$.
\end{lem}

\begin{proof}
  By assumption, we have that $L|_{\ell_i}\simeq \ring \P^1.(d_i)$.
  Let $\inter$ be the scheme theoretic intersection of $\ell_1$ and
  $\ell_2$, $\iota^i:\ell_i\to\ell$ the natural embedding, and
  $L_i=\iota^i_*(L|_{\ell_i})$ for $i=1,2$.  Then one has the
  following short exact sequence:
  $$
  0\to L \to L_1\oplus L_2 \to \ring \inter. \to 0.
  $$
  This implies that $h^0(\ell, L)\geq
  \chi(L)=\chi(L_1)+\chi(L_2)-\chi(\ring\inter.)= d_1+d_2$.
\end{proof}

\subsection{The identification of the components of a dubby}

To illustrate the main observation about dubbies, let us consider a
very simple setup first: let $L \in \Pic(X)$ be an ample line bundle,
and assume that $\ell = \ell_1 \cup \ell_2 \subset X$ is a dubby where
both components are members of the same connected family $H$ of
minimal rational curves. In particular, $L|_\ell$ will be of type
$(d,d)$, where $d>0$. Remarkably, the line bundle $L$ induces a
canonical identification of the two components $\ell_1$ and $\ell_2$,
at least when $d$ is coprime to the characteristic of the base field
$k$. Over the field of complex numbers, the idea of construction is
the following: Fix a trivialization $t:L|_V\to \O_V$ of $L$ on an open
neighborhood $V$ of the intersection point $\{z\} = \ell_1 \cap\,
\ell_2$.  Given a point $x \in \ell_1 \setminus \ell_2$, let $\sigma_1
\in H^0(\ell_1, L|_{\ell_1})$ be a non-zero section that vanishes at
$x$ with multiplicity $m$.  Then there exists a \emph{unique} section
$\sigma_2 \in H^0(\ell_2, L|_{\ell_2})$ with the following properties:
\begin{enumerate} 
\item The section $\sigma_2$ vanishes at exactly one point $y \in
  \ell_2$.
\item The sections $\sigma_1$ and $\sigma_2$ agree on the intersection
  of the components:
  $$
  \sigma_1(z) = \sigma_2(z)
  $$
\item The differentials of $\sigma_1$ and $\sigma_2$ agree at $z$:
  $$
  \vec v(t\circ \sigma_1) = \vec v(t\circ \sigma_2)
  $$
for all non-vanishing tangent vectors $\vec v\in T_{\ell_1}\cap
T_{\ell_2} $.
\end{enumerate}
The map that associates $x$ to $y$ gives the identification of the
components and does not depend on the choice of $t$.

In the following section~\ref{sec:bundles_of_dubbies}, we will give a
construction of the identification morphism which also works in the
relative setup, for bundles of type $(d_1, d_2)$ where $d_1 \ne d_2$,
and in arbitrary characteristic.

\subsection{Bundles of dubbies}
\label{sec:bundles_of_dubbies}

For the proof of the main theorems we will need to consider bundles of
dubbies, i.e., morphisms where each scheme-theoretic fiber is
isomorphic to a dubby. The following Proposition shows how to identify
the components of such bundles.

\begin{prop}\label{prop:identification_of_dubbies-epsilon}
  Let $\lambda: \Lambda \to B$ be a projective family of dubbies over
  a normal base $B$ and assume that $\Lambda$ is not irreducible.
  Then it has exactly two irreducible components $\Lambda_1$ and
  $\Lambda_2$, both $\P^1$-bundles over $B$. Assume further that there
  exists a line bundle $L \in \Pic(\Lambda)$ whose restriction to a
  $\lambda$-fiber has type $(m, n')$, where $m$ and $n$ are non-zero and
  relatively prime to $\charact(k)$.
  
  If $\Sigma_{\red} \subset \Lambda_1 \cap \Lambda_2$ denotes the
  reduced intersection, then $\Sigma_{\red}$ is a section over $B$, and
  the pairs $(\Lambda_1, \Sigma_{\red})$ and $(\Lambda_2,
  \Sigma_{\red})$ are isomorphic over $B$.
\end{prop}

Note that the isomorphism given in
Proposition~\ref{prop:identification_of_dubbies-epsilon} is not
canonical and may not respect the line bundle $L$.

 \subsubsection*{Proof of
Proposition~\ref{prop:identification_of_dubbies-epsilon}}

The map $\lambda$ is flat because all its scheme-theoretic fibers are
isomorphic. Let $\Lambda_1 \subset \Lambda$ be one of the irreducible
components. It is easy to see that if $x\in\Lambda_1$ is a general
point, then $\Lambda_1$ contains the (unique) irreducible component of
$\ell_{\lambda(x)} := \lambda^{-1} \lambda(x)$ that contains $x$.
Since $\lambda$ is proper and flat, $\lambda(\Lambda_1)=B$. Hence
$\Lambda_1$ contains one of the irreducible components of $\ell_b$ for
all $b\in B$.  Repeating the same argument with another irreducible
component, $\Lambda_2$, one finds that it also contains one of the
irreducible components of $\ell_b$ for all $b\in B$. However, they
cannot contain the same irreducible component for any $b\in B$: In
fact, if they contained the same component of $\ell_b$ for infinitely
many points $b\in B$, then they would agree. On the other hand, if
they contained the same component of $\ell_b$ for finitely many points
$b\in B$, then $\Lambda$ would have an irreducible component that does
not dominate $B$. This, however, would contradict the flatness of
$\lambda$.  Hence $\Lambda_1 \cup \Lambda_2 = \Lambda$. They are both
$\P^1$-bundles over $B$ by \cite[Thm.~II.2.8.1]{K96}.

Let $\Sigma := \Lambda_1 \cap \Lambda_2$ be the scheme-theoretic
intersection. Since $\Lambda$ is a bundle of dubbies and $B$ is
normal, it is clear that its reduction, $\Sigma_{\red}$ is a section,
and that $\Sigma$ is its first infinitesimal neighborhood in either
$\Lambda_1$ or $\Lambda_2$.  In this setup, the isomorphism of pairs
is given by Theorem~\ref{thm:isomorphism_of_pairs}. \qed

\subsection{The space of dubbies}
\label{sect:space}

 In addition to the space of rational curves, which we use throughout,
it is also useful to have a parameter space for dubbies.  For the
convenience of the reader, we will first recall the construction of the
former space very briefly.  The reader is referred to
\cite[chapt.~II.1]{K96} for a thorough treatment.

\subsubsection{The space of rational curves}
\label{sec:ratcurves-review}

Recall that there exists a scheme $\Hom_{\operatorname{bir}} (\P^1,
X)$ whose geometric points correspond to morphisms $\P^1 \to X$ that
are birational onto their images. Furthermore, there exists an
``evaluation morphism'': $\mu: \Hom_{\operatorname{bir}}(\P^1,X)\times
\P^1 \to X$. The group $\psl_2$ acts on the normalization
$\Hombir(\P^1, X)$, and the geometric quotient exists. More precisely,
we have a commutative diagram
$$
  \label{diag:rat_curves}
  \xymatrix{ 
    \Hombir (\P^1, X)\times \P^1 \ar[d] \ar[r]^(.6){U} 
    \ar@/^.6cm/[rr]^{\mu}
    & {\Univrc(X)} \ar[r]^(.6){\iota} \ar[d]_{\pi} & {X} \\ 
    \Hombir (\P^1, X) \ar[r]^{u} & {\RatCurves(X)} }  
$$
where $u$ and $U$ are principal $\psl_2$ bundles, $\pi$ is a
$\P^1$-bundle and the restriction of the ``evaluation morphism''
$\iota$ to any fiber of $\pi$ is a morphism which is birational onto
its image. The quotient space $\RatCurves(X)$ is then the parameter
space of rational curves on $X$.  The letter ``$\operatorname{n}$'' in
$\RatCurves$ may be a little confusing. It has nothing to do with the
dimension of $X$ and it's not a power. It serves as a reminder that
the parameter space is the $\operatorname{n}$ormalization of a
suitable quasiprojective subset of the Chow variety.

It may perhaps look tempting to define a space of dubbies in a similar
manner, as a quotient of the associated Hom-scheme. However, since
geometric invariant theory becomes somewhat awkward for group actions
on non-normal varieties, we have chosen another, elementary but
somewhat lengthier approach. The space of dubbies will be constructed
as a quasi-projective subvariety of the space of ordered pairs of
pointed rational curves, and the universal family of dubbies will be
constructed directly.

\subsubsection{Pointed rational curves}
It is easy to see that $\ptRC(X) =\Univrc(X)$ naturally parametrizes
pointed rational curves on $X$ and the pull-back of the universal
family 
$$
\ptUnivrc(X)=\ptRC(X)\times_{\RC(X)}\Univrc(X)
$$ is the universal family of pointed rational curves over
$\ptRC(X)$. The identification morphism $\ptRC(X)\to\Univrc(X)$ and
the identity map of $\ptRC(X)$ gives a section of this universal
family:
$$
\xymatrix{ 
{\qquad\ptUnivrc(X)} \ar[d]  \ar[r] & \Univrc(X) \ar[d] \\ 
{\qquad\ptRC(X)} \ar@{..>}@/^0.5cm/[u]^{\eta} 
\ar[ru]_{\simeq} \ar[r]_\pi & {\RC(X)}
}
$$

\subsubsection{Pairs of pointed rational curves}

The product $\twoptRC(X) := \ptRC(X)\times\ptRC(X)$ naturally
parametrizes pairs of pointed rational curves. We denote the
projections to the two factors by $\rho_i:\twoptRC(X)\to\ptRC(X)$ for
$i=1,2$. Then the universal family will be given as the disjoint union
$$
\twoptUnivrc(X)=\left(\twoptRC(X)
  \times_{\rho_1}\ptUnivrc(X)\right)\cup
\left(\twoptRC(X)\times_{\rho_2}\ptUnivrc(X)\right).
$$
The two copies of the section $\eta:\ptRC(X)\to \ptUnivrc(X)$
induce two sections of this family, one for each component of the
union:
$$
  \label{diag:2ptrat_curves} \xymatrix{ & \ar[dl] {\twoptUnivrc(X)}
 \ar[dd]_{\tilde p} \ar [dr] & \\
\ptUnivrc(X) \ar[dd] && \ptUnivrc(X) \ar [dd] \\
 & \ar[dl]_{\rho_1}  \twoptRC(X)\ar [dr]^{\rho_2} \
\ar@{..>}@/^0.5cm/[uu]^{\sigma_1} 
\ar@{..>}@/_0.5cm/[uu]_{\sigma_2} 
&\\
\ar@{..>}@/^0.5cm/[uu]^{\eta} 
\ptRC(X) && 
\ar@{..>}@/_0.5cm/[uu]_{\eta}
\ptRC(X)
  }  
$$

\subsubsection{The space of dubbies}
\label{sec:defn_of_dubbies}

Consider the evaluation morphism $\iota^2: \twoptUnivrc(X) \to X$. The
associated tangent map $T\iota^2$ restricted to the relative tangent
sheaf $T_{\twoptUnivrc(X)/\twoptRC(X)}$ gives rise to a rational map
$$
\twopttau
 : \twoptUnivrc(X) \dasharrow \P(T_X^\vee).
$$
We define a quasiprojective variety, the space of dubbies,
\begin{multline*}
  \Dubbies(X) := \text{normalization of }\{\ell \in \twoptRC(X) \,|\, \text{$\twopttau$ is defined at $\sigma_1(\ell)$} \\
  \text{and at $\sigma_2(\ell)$, and $\twopttau(\sigma_1(\ell))=\twopttau(\sigma_2(\ell))$}\}.
\end{multline*}
We will often consider pairs of curves such that both components come
from the same family $H \subset \RatCurves(X)$. For this reason we
define $\pi_2:\Dubbies(X) \to \RatCurves(X)\times\RatCurves(X)$, the
natural forgetful projection morphism, and
$$
\Dubbies(X)|_H := \Dubbies(X) \cap \pi_2^{-1}(H\times H).
$$

\begin{prop}\label{prop:properness-of-dubbies}
  Assume that $H \subset \RatCurves(X)$ is a proper family of
  immersed curves. Then $\Dubbies(X)|_H$ is also proper.
\end{prop}

\begin{proof}
Since the tangent map, $\twopttau$, is well-defined at
  $\sigma_1(\ell)$ and $\sigma_2(\ell)$ for every $\ell \in
  \twoptRC(X) \cap\, \pi_2^{-1}(H\times H)$,
  $$
  \Dubbies(X) = \text{normalization of }\{\ell \in \twoptRC(X) \,|\,
  \twopttau(\sigma_1(\ell)) = \twopttau(\sigma_2(\ell)) \},
  $$
  which is clearly a closed subvariety of the proper variety
  $\pi_2^{-1}(H\times H)$.
\end{proof}

The next statement follows immediately from the construction and from
the universal property of $\RatCurves(X)$.

\begin{prop}\label{prop:univProp}
  Let $\ell_1$ and $\ell_2 \subset X$ be rational curves with
  normalizations 
  $$
  \eta_i:\P^1\simeq \tilde\ell_i \to \ell_i\subset X.
  $$
  If $T\eta_i$ have rank 1 at the point $[0:1]\in \P^1$ for $i=1,2$,
  and if the images of the tangent morphisms agree,
  $$
  \Image(T\eta_1|_{[0:1]}) = \Image(T\eta_2|_{[0:1]}) \subset T_X,
  $$
  then there exists a point $\ell \in \Dubbies(X)$ such that $\tilde
  p^{-1}(\ell) = \tilde\ell_1\cup\tilde\ell_2$.
  
  If $H \subset \RatCurves(X)$ is a subfamily, and both $\ell_i$
  correspond to points of $H$, then we can find such an $\ell$ in
  $\Dubbies(X)|_H$. \qed
\end{prop}

\subsubsection{The universal family of dubbies}

In order to show that $\Dubbies(X)$ is a space of dubbies indeed, we
need to construct a univeral family, which is a bundle of dubbies in
the sense of section~\ref{sec:bundles_of_dubbies}. To this end, we
will factor the universal evaluation morphism.

\begin{prop}\label{prop:triviality_of_U}
  The evaluation morphism,
  $$ \iota :
  \underbrace{\twoptUnivrc(X)\times_{\twoptRC(X)}{\Dubbies(X)}}_{=:
  \tilde U, \text{ decomposes as } \tilde U_1 \cup \tilde U_2}
  \longrightarrow U \subset X\times \Dubbies(X),
  $$
  factors as follows:
  \begin{equation}
    \label{eq:family_of_dubbies}
  \xymatrix{ {\tilde U} \ar[rr]^{\alpha}
    \ar@/^0.8cm/[rrrr]^{\iota} \ar@/_/[drr]_{\txt{\scriptsize two disjoint \quad \quad \\
        \scriptsize $\P^1$-bundles $\tilde U_1 \cup \tilde U_2$\quad
        \quad}}^{\tilde p} & & {\Lambda}
    \ar[rr]^{\beta} \ar[d]^{\hat p}_{\txt{\scriptsize bundle of \\
        \scriptsize dubbies}} & & {U}
    \ar@/^/[dll]_{p}^{\txt{\scriptsize \quad \quad \quad bundle of two curves with \\
        \scriptsize \quad \quad \quad complicated intersection}} \\ &
    & {\Dubbies(X)} & }
  \end{equation}
  where the variety $\Lambda$ is a bundle of dubbies in the sense
  that for every closed point $b\in B^0$, the fiber $\hat p^{-1}(b)$
  is isomorphic to a dubby.
\end{prop}

\begin{rem}
  If $\ell \in \Dubbies(X)$ is any point, then the two corresponding
  curves in $X$ intersect tangentially in one point, but may have very
  complicated intersection at that point and elsewhere. The
  factorization of the evaluation morphism should therefore be
  understood as a partial resolution of singularities.
\end{rem}

\begin{proof}
  As a first step we will construct the space $\Lambda$. Because the
  evaluation $\iota$ is a finite, hence affine, morphism, it seems
  appropriate to construct a suitable subsheaf $\mathcal A \subset
  \iota_*\mathcal O_{\tilde U}$, which is a coherent sheaf of
  $\mathcal O_U$-modules and set $\beta: \Lambda =\sheafspec(\mathcal
  A)\to U$.
  
  Let $\bar\sigma_1 \subset \tilde U_1$ and $\bar\sigma_2 \subset
  \tilde U_2$ be the images of the pullbacks of the canonical
  sections, $\sigma_1$ and $\sigma_2$, constructed in
  \ref{diag:2ptrat_curves}. In order to construct $\mathcal A$, we
  will need to find an identification of their first infinitesimal
  neighborhoods, $\tilde \sigma_1$ and $\tilde \sigma_2$. Since
  $\iota$ is separable, it follows directly from the construction that
  $\tilde \sigma_1$ and $\tilde \sigma_2$ map isomorphically onto
  their scheme-theoretic images $\iota(\tilde \sigma_1)$ and
  $\iota(\tilde \sigma_2)$. Again, by the definition of $\Dubbies(X)$,
  these images agree: $\iota(\tilde \sigma_1) = \iota(\tilde
  \sigma_2)$ and we obtain the desired identification,
  $$
  \gamma : \tilde \sigma_1 \to \tilde \sigma_2.
  $$
  Let 
  $$
  i_1 : \tilde \sigma_1 \to \tilde U_1 \quad \text{and} \quad i_2 :
  \tilde \sigma_2 \to \tilde U_2
  $$
  be the inclusion maps and consider the sheaf morphism
  $$
  \varphi := \iota_*(i_1^\#-\gamma^\# \circ i_2^\#) : \iota_*
    \O_{\tilde U} \to \iota_*\O_{\tilde \sigma_1}.
  $$
  The sheaf 
  $$
  \mathcal A := \ker(\varphi)
  $$ is thus a coherent sheaf of $\mathcal O_U$-modules. As it was
  planned above, define $\Lambda := \sheafspec (\mathcal A)$. The
  existence of the morphisms $\alpha$ and $\beta$ and that
  $\iota=\beta\circ\alpha$ follows from the construction. It remains
  to show that $\Lambda$ is a bundle of dubbies.  Let $\ell \in
  \Dubbies(X)$ be a closed point. Replacing $\Dubbies(X)$ with a
  neighbourhood of $\ell$ and passing to a finite, unbranched cover if
  necessary, and by abuse of notation still denoting it by
  $\Dubbies(X)$, we can assume that
  \begin{enumerate} 
  \item the variety $\Dubbies(X)$ is affine, say $\Dubbies(X) \simeq
    \Spec \ringB$,
  \item the $\P^1$-bundles $\tilde U_i = \P(\tilde p_* \O_{\tilde
      U_i}(\bar\sigma_i)^\vee)$, for $i=1,2$ are trivial, and
  \item there exists a Cartier divisor $\tau \subset U$ such that
    $\iota^{-1}(\tau) = \tau_1 \cup \tau_2$, where $\tau_i \subset
    \tilde U_i$ are sections that are disjoint from $\bar\sigma_i$.
  \end{enumerate} 
  We can then find homogeneous bundle coordinates $[x_0:x_1]$ on
  $\tilde U_1$ and $[y_0:y_1]$ on $\tilde U_2$ such that
  \begin{align*} \bar\sigma_1 & = \{ ([x_0:x_1],b) \in \tilde U_1
    \,|\, x_0 =0 \}, &
    \tau_1 & = \{ ([x_0:x_1],b) \in \tilde U_1 \,|\, x_1 =0 \}, \\
    \bar\sigma_2 & = \{ ([y_0:y_1],b) \in \tilde U_2 \,|\, y_0 =0 \},
    \text{ and} & \tau_2 & = \{ ([y_0:y_1],b) \in \tilde U_2 \,|\, y_1
    =0 \}.  \end{align*} If we set
  $$
  \tilde U_0 := \tilde U \setminus (\tau_1 \cup \tau_2),
  $$
  then the image $U_0 := \iota(\tilde U_0)$ is affine, and we can
  write the relevant modules as
  \begin{align*} 
    \O_{\tilde U}(\tilde U_0) & \simeq \ringB \otimes (k[x_0] \oplus
    k[y_0]) \\
    \O_{\tilde \sigma_1}(\tilde U_0) & \simeq \ringB \otimes 
    k[x_0] / (x_0^2).\\
    \O_{\tilde \sigma_2}(\tilde U_0) & \simeq \ringB \otimes 
    k[y_0] / (y_0^2).
  \end{align*} 
  Adjusting the bundle coordinates, if necessary, we can assume that
  the identification morphism $\gamma^\#(U_0) : \O_{\tilde
    \sigma_2}(\tilde U_0) \to \O_{\tilde \sigma_1}(\tilde U_0)$ is
  written as
  $$
  \begin{array}{rccc}
    \gamma^\#(U_0) : & \ringB \otimes k[y_0]/(y_0^2) & \to & \ringB
    \otimes k[x_0] / (x_0^2) \\ & r \otimes y_0 & \mapsto & r
    \otimes x_0.
  \end{array}
  $$
  In this setup, we can find the morphism $\varphi$ explicitly:
  $$
  \begin{array}{rccc} 
    \varphi(U_0) : & \ringB \otimes (k[x_0] \oplus
    k[y_0]) & \to & \ringB \otimes k[x_0] / (x_0^2) \\ 
    & r\otimes (f(x_0),g(y_0)) & \mapsto & r \otimes (f(x_0)-g(x_0)).
  \end{array}
  $$
  Therefore, as an $\ringB$-algebra, $\ker(\varphi)(U_0)$ is
  generated by the elements $u := 1_{\ringB} \otimes (x_0,y_0)$ and $v
  := 1_{\ringB} \otimes (x_0^2,0)$, which satisfy the single relation
  $v(u^2-v) =0$. Thus
  $$
  \ker (\varphi)(U_0) = \ringB \otimes k[u,v] / (v(u^2-v)).
  $$
  In other words, $\beta^{-1}(U_0)$ is a bundle of two affine lines
  over $\Dubbies(X)$,
  meeting tangentially in a single point.
  
  It follows directly from the construction of $\mathcal A$ that
  $\alpha$ is an isomorphism away from $\bar\sigma_1 \cup \bar\sigma_2$.
  The curve $\hat p^{-1}(\ell)$ is therefore smooth outside of $\hat
  p^{-1}(\ell) \cap \beta^{-1}(U_0)$, and it follows that $\hat
  p^{-1}(\ell)$ is indeed a dubby. This ends the proof.
\end{proof}

\section{Proofs of the Main Theorems}

\subsection{Proof of Theorem~\ref{thm:main-criterion1}}
\label{chap:proof11}

The assertion that all curves associated with $H_x$ are smooth at a
general point $x \in X$ follows immediately from the assumption that
none of the curves $\ell \in H$ is cuspidal, and by \cite[thms.~2.4(1)
and 3.3(1)]{Keb01b}. It remains to show that no two curves intersect
tangentially.

We will argue by contradiction and assume that we can find a pair
$\ell = \ell_1 \cup \ell_2 \subset X$ of distinct curves $\ell_i \in
H$ that intersect tangentially at $x$.  The pair $\ell$ is then
dominated by a dubby whose singular point maps to $x$. Loosely
speaking, we will move the point of intersection to obtain a
positive-dimensional family of dubbies that all contain the point $x$.

\subsubsection*{Setup}

To formulate our setup more precisely, we will use the notation
introduced in diagram~(\ref{eq:family_of_dubbies}) of
Proposition~\ref{prop:triviality_of_U} and recall from
Proposition~\ref{prop:properness-of-dubbies} that $\Dubbies(X)|_H$ is
proper. Recall further that the universal family $U$ is a subset $U
\subset X \times \Dubbies(X)$ and let $\jmath := pr_1\circ \iota : \tilde U
\to X$ be the canonical morphism. The assumption that for every
general point $x \in X$, there is a pair of curves intersecting
tangentially at $x$ can be reformulated as
$$
\jmath \circ \sigma_1(\Dubbies(X)|_H) = \jmath \circ \sigma_2(\Dubbies(X)|_H) = X
$$
Let $D \subset \Dubbies(X)|_H$ be an irreducible component such that
$$
\jmath \circ \sigma_1(D) = \jmath \circ \sigma_2(D) = X
$$ holds.  By abuse of notation, we will denote $\tilde U_D=(\tilde
U_{D})_1\cup (\tilde U_{2})_D$ by $\tilde U=\tilde U_1\cup \tilde
U_2$. Fix a closed point $t \in D$ and consider the intersection
numbers
$$
d_1 := \jmath^*(L)\cdot (\tilde p^{-1}(t)\cap \tilde U_1) 
\quad \text{and} \quad
d_2 := \jmath^*(L)\cdot (\tilde p^{-1}(t)\cap \tilde U_2) 
$$
Renumbering $\tilde U_1$ and $\tilde U_2$, if necessary, we may
assume without loss of generality that $d_1 \geq d_2$. In this setup
it follows from the upper semi-continuity of the fiber dimension that
$(\jmath|_{\tilde U_1})^{-1}(x)$ contains an irreducible curve $\tau_1$
which intersects $\sigma_1(D)$ non-trivially and is not contained in
$$
S := \{ y \in \tilde U \,|\,  \text{$\iota$ is not an isomorphism at $y$}\}
$$
Set $T:=\tilde{p}(\tau_1)$. After a base change, if necessary, we
may assume that $T$ is a normal curve and consider the restrictions of
the morphisms constructed in Proposition~\ref{prop:triviality_of_U}:
$$
\xymatrix{ 
  {\tilde U_T} \ar[r]_{\alpha} \ar@/^0.5cm/[rrr]^{\jmath}
  & {\Lambda_T} \ar[r]_{\beta} & {U_T} \ar[r]_{pr_1} & X
}
$$ Using \cite[thm.~3.3.(1)]{Keb00a}, we find that $\tau_1$ is
generically injective over $T$, and therefore is a section.  Let
$\tilde U_{T,1}=(\tilde U_1)_T$ and $\tilde U_{T,2}=(\tilde
U_2)_T$. By Proposition~\ref{prop:identification_of_dubbies-epsilon},
$(\tilde U_{T,1}, \sigma_1(T))$ and $(\tilde U_{T,2}, \sigma_2(T))$
are isomorphic pairs over $T$. Let $\gamma :\tilde U_{T,1} \to \tilde
U_{T,2}$ be an isomorphism and consider the section $\tau_2 :=
\gamma(\tau_1) \subset \tilde U_{T,2}$.

\subsubsection*{The contraction of $\tau_2$}

With the notation above, Theorem~\ref{thm:main-criterion1} follows
almost immediately from the following observation.

\begin{lem}\label{lem:contr_of_tau2}
  The morphism $\jmath$ contracts the section $\tau_2$ to $x$, i.e.,
  $\tau_2 \subset \jmath^{-1}(x)$.
\end{lem}

Notice that this finishes the proof of
Theorem~\ref{thm:main-criterion1}. Indeed,
Lemma~\ref{lem:contr_of_tau2} implies that a general point $t\in T$
corresponds to a pair $\ell_t = \ell_{t,1}\cup \ell_{t,2}$ of two
distinct curves that intersect at $x$. The curve $\ell_t$ is then
singular at $x$, a contradiction to the fact that $\tau_1 \not \subset
S$.

\begin{proof}[Proof Lemma~\ref{lem:contr_of_tau2}]
  As a first step, we show that $\jmath$ contracts $\tau_2$ to some
  point $y \in X$. The proof relies on a calculation of intersection
  numbers on the ruled surfaces $\tilde U_{T,1}$ and $\tilde U_{T,2}$.
  Recall the basic fact that
  $$
  \Num(\tilde U_{T,1}) =\ZZ\cdot \sigma_1(T) \oplus \ZZ \cdot F_{V,1}
  $$ where $F_{V,1}$ is a fiber of $\tilde p_{\tilde U_{T,1}}:\tilde
  U_{T,1} \to T$. A similar decomposition holds for $\tilde
  U_{T,2}$. Since $\tau_1$ is a section, we have the numerical
  equivalence,
  $$
  \tau_1 \equiv \sigma_1(T) + d\cdot F_{V,1},
  $$
  where $d$ is a suitable integer. Since $\gamma$ maps
  $\sigma_1(T)$ isomorphically onto $\sigma_2(T)$, we obtain a similar
  equation on $\tilde U_{T,2}$,
  $$
  \tau_2 \equiv \sigma_2(T) + d\cdot F_{V,2}.
  $$ 
Next take the ample line bundle $L\in \Pic(X)$ and set
  $$
  d_3 := \jmath^*(L)\cdot\sigma_1(T) =\jmath^*(L)\cdot\sigma_2(T).
  $$ These two numbers are indeed equal since the evaluation morphism
identifies the images of the two sections $\sigma_1(T)$ and
$\sigma_2(T)$. Now we can write the intersection numbers as
  $$
  \begin{aligned}
    \jmath^*(L) \cdot \tau_2 & = \jmath^*(L) \cdot( \sigma_2(T) + d\cdot F_{V,2} ) \\
    & = d_3+d \cdot d_2 = \underbrace{(d_3 + d\cdot d_1)}_{= \jmath^*(L)
      \cdot \tau_1 = 0} + d\cdot (d_2-d_1) \\
    & = d\cdot (d_2-d_1) \leq 0
  \end{aligned}
  $$ Since $L$ is ample, this shows that $\jmath(\tau_2)$ is a point,
  $y \in X$.
  
  It remains to prove that $x=y$. In order to see that, it suffices to
  recall two facts. First, as it was already used above, the
  evaluation morphism identifies the images of the two sections
  $\sigma_1(T)$ and $\sigma_2(T)$. Second, we know that $\tau_1$ and
  the canonical section $\sigma_1 \subset \Lambda_1$ intersect. Let $t
  \in \tilde p(\tau_1 \cap \sigma_1(T))$ be a closed point. The two
  sections $\tau_2$ and $\sigma_2(T)$ will then also intersect, $t \in
  \tilde p(\tau_2 \cap \sigma_2(T))$ and we obtain
  \begin{align*}
    x & = \jmath(\tau_1) = \jmath(\tau_1 \cap \sigma_1(T) \cap \tilde p^{-1}(t)) \\
    & = \jmath(\tau_2 \cap \sigma_2(T) \cap \tilde p^{-1}(t)) = y.
  \end{align*}
  Lemma~\ref{lem:contr_of_tau2} is thus shown.
\end{proof}

\subsection{Proof of Theorem~\ref{thm:main-criterion2}}
\label{chap:proof12}

Let $H \subset \RatCurves(X)$ be as in
Theorem~\ref{thm:main-criterion2}. We assume without loss of
generality that all irreducible components of $H$ dominate $X$. Fix an
ample line bundle $L \in \Pic(X)$ and let $H' \subset H$ be an
irreducible component such that for a general curve $\mathcal C \in
H'$ the intersection number $L.\mathcal C$ is minimal among all the
intersection numbers of $L$ with curves from $H$.  Finally, fix a
rational curve $\mathcal C \subset X$ that corresponds to a general
point of $H'$.

The proof of Theorem~\ref{thm:main-criterion2} now follows very much
the lines of the proof of Theorem~\ref{thm:main-criterion1} from the
previous section.  The main difference to the previous argument is
that we have to work harder to find the family $T$, as the properness
of $\Dubbies(X)|_H$ is no longer automatically guaranteed. Over the
complex number field, however, the following lemma holds, which
replaces the properness assumption in our context.
\begin{lem}
  Assume that $X$ is a complex-projective manifold, and let $S'\subset
  X$ be a subvariety of codimension $\codim_XS'\geq 2$.  If $\mathcal
  C \in H$ is a curve that corresponds to a general point of $H'$,
  then $\mathcal C$ and $S'$ are disjoint: $\mathcal C \cap
  S'=\emptyset$.
\end{lem}
\begin{proof}
  \cite[Chapt.~II, Prop.~3.7 and Thm.~3.11]{K96}
\end{proof}

\begin{cor}\label{cor:properness_of_immersed}
  Under the assumptions of Theorem~\ref{thm:main-criterion2}, if
  $\mathcal C\in H'$ is a general curve, and if $\codim_XD\geq 2$,
  then
  $$
  H_{\mathcal C} := \{ \mathcal C' \in H \,|\, \mathcal
  C\cap\mathcal C'\neq\emptyset\} \subset H
  $$
  is proper, and the associated curves are immersed. In particular,
  $\mathcal C$ is immersed.
\end{cor}
\begin{proof}
  It suffices to note that $\mathcal C$ is disjoint from both $S$ and
  $D$.
\end{proof}

Before coming to the proof of Theorem~\ref{thm:main-criterion2}, we
give a last preparatory lemma concerning the dimension of the locus
$D$ of cusps.

\begin{lem}\label{lem:dimension_of_D}
  If $D \subset X$ is a divisor, then the subfamily $H^{\text{cusp}}
  \subset H$ of cuspidal curves dominates $X$.
\end{lem}
\begin{proof}
  Argue by contradiction and assume that all cuspidal curves in
  $H^{\text{cusp}}$ are contained in a divisor. The total space of the
  family of cuspidal curves is then at least $(\dim D+1)$-dimensional,
  so for a general point $x\in D$ there exists a positive dimensional
  family of cuspidal curves that contain $x$ and are contained in $D$.
  That, however, is impossible: it has been shown in
  \cite[Thm.~3.3]{Keb00a} that in the projective variety $D$, a
  general point is contained in no more then finitely many cuspidal
  curves.
\end{proof}

\subsubsection*{Setup of the proof}

For the proof of Theorem~\ref{thm:main-criterion2}, we will again
argue by contradiction. By Lemma~\ref{lem:dimension_of_D} this amounts
to the assumtion that $\tau$ is not generically injective, and that
$\codim_X D\geq 2$. By Corollary~\ref{cor:properness_of_immersed},
this implies that the space of curves which intersect $\mathcal C$ is
proper and all associated curves are immersed along $\mathcal C$.
Since $\mathcal C$ was a general curve, the assumptions also imply
that for a general point $x \in \mathcal C$, there exists a point $t
\in \Dubbies(X)$ corresponding to a pair of marked curves $\ell =
\ell_1 \cap \ell_2$ such that $\ell_2 = \mathcal C$ and $\ell_1$
intersects $\mathcal C$ tangentially at $x$, i.e.,
$\Image(\tau(\sigma_1(t)))= \P(T_{\mathcal C}|_x^\vee)$ where
$\tau:\tilde U\to \P(T_X^{\vee})$ is the tangent morphism from the
introduction.

Hence we can find a proper curve $T\subset \Dubbies(X)$ with
associated diagram
$$
\xymatrix{ {\tilde U_T} \ar[r]_{\alpha} \ar@/^0.5cm/[rrr]^{\jmath}
  \ar@/_/[dr]_{\tilde p} & {\Lambda_T} \ar[r]_{\beta} \ar[d]^{\hat p}
  & {U_T} \ar[r]_{pr_1}
  \ar@/^/[dl]^{p} & X\\
  & {T} & }
$$ 
such that $\tilde U_T$ decomposes as $\tilde U_T=\tilde U_{T,1}
\cup \tilde U_{T,2}$, where
$$
\tilde U_{T,2} \simeq \tilde {\mathcal C} \times T \simeq \P^1 \times
T,
$$
and where $\tau|_{\sigma_1(T)}$ dominates $\P(T_{\mathcal C}^\vee)$.

\subsubsection*{End of proof}

We are now in a situation which is very similar to the one considered
in the proof of Theorem~\ref{thm:main-criterion1}: we will derive a
contradiction by calculating certain intersection numbers on $\tilde
U_{T,1}$ and $\tilde U_{T,2}$.

As a first step, remark that $\tilde U_{T,1}$ maps to a surface in
$X$. It follows that $\jmath^*(L)$ is nef and big on $\tilde U_{T,1}$.

Secondly, since $\tilde U_{T,2}$ is isomorphic to the trivial bundle $\P^1
\times T$, we have a decomposition
$$
\Num(\tilde U_{T,2}) \simeq \ZZ \cdot F_{H,2} \oplus \ZZ \cdot F_{V,2}
$$
where $F_{H,2}$ is the numerical class of a fiber of the map
$\tilde U_{T,2} \to \P^1$ and $F_{V,2}$ that of a fiber of the map
$\tilde U_{T,2} \to T$. Likewise, since the pairs $(\tilde
U_{T,1}, \sigma_1(T))$ and $(\tilde U_{T,2}, \sigma_2(T))$ are
isomorphic, let 
$$
\Num(\tilde U_{T,1}) \simeq \ZZ \cdot F_{H,1} \oplus \ZZ \cdot F_{V,1}
$$
be the corresponding decomposition. If $d$ denotes the degree of the
(finite, surjective) morphism
$$
\jmath \circ \sigma_1 = \jmath \circ \sigma_2 : T \to \mathcal C,
$$
then it follows directly from the construction that the curves of type
$F_{H,2}$ intersect $\sigma_2(T)$ with multiplicity $d$. We obtain
that
$$
\sigma_2(T) \equiv F_{H,2}+d\cdot F_{V,2}  \quad \text{and thus} \quad
\sigma_1(T) \equiv F_{H,1}+d\cdot F_{V,1}.
$$
To end the argumentation, let
$$
d_1 := \jmath^*(L)\cdot F_{V,1} \quad \text{and} \quad d_2 :=
\jmath^*(L)\cdot F_{V,2}
$$
In particular, we have that $\jmath^*(L)\cdot \sigma_2 = d \cdot k_2$.
Recall that $H' \subset H$ was chosen so that $d_1 \geq d_2$ and
write:
$$
\begin{aligned}
  \jmath^*(L)\cdot F_{H,1} & = \jmath^*(L)\cdot (\sigma_1(T) - d\cdot F_{V,1}) \\
  & = d \cdot d_2 - d \cdot d_1 \\
  & \leq 0.
\end{aligned}
$$
Because $\tilde U_{T,1}$ is covered by curves which are numerically
equivalent to $F_{H,1}$ that contradicts the assumption that
$\jmath^*(L)|_{\tilde U_{T,1}}$ is big and nef.  The proof of
Theorem~\ref{thm:main-criterion2} is thus finished. \qed

\section{Applications}\label{sect:applications}

\subsection{Irreducibility Questions}

Let $H \subset \RatCurves(X)$ be a maximal dominating family of
rational curves of minimal degrees on a projective variety $X$. How
many components can $H$ have? If we pick an irreducible component $H'
\subset H$ and fix a general point $x\in X$, does it follow that
$$
H'_x := \{\ell \in H'\,|\, x\in \ell\}
$$ is irreducible? These questions have haunted the field for quite a
while now, as the possibility that $H'_x$ might be reducible poses
major problems in many of the proposed applications of rational curves
to complex geometry ---see the discussion in \cite{Hwa}.

It is conjectured \cite[chap.~5, question 2]{Hwa} that the answers to
both of these questions are affirmative for a large class of
varieties. There exists particularly strong evidence if $X$ is a
complex manifold and if the dimension of $H'_x$ is not too small.
Theorem~\ref{thm:main-criterion1} enables us to give a partial answer.

\begin{thm}\label{thm:irred}
  Under the assumptions of Theorem~\ref{thm:main-criterion1}, if $X$
  is a complex manifold and if for a general point $x\in X$, and for
  all irreducible components $H' \subset H$
  $$
  \dim H'_x \geq \frac{\dim X-1}{2},
  $$
  then $H_x$ is irreducible. In particular, $H$ is irreducible.
\end{thm}

The main technical difficulty in proving Theorem~\ref{thm:irred} lies
in the fact that the closed points of $H$ are generally not in
1:1-correspondence with actual rational curves, a possibility that is
sometimes overlooked in the literature. As a matter of fact, this
correspondence is only generically injective, and it may well happen
that two or more points of $H$ correspond to the same curve $\ell
\subset X$. This is due to the very construction of the space
$\RatCurves(X)$: recall from section~\ref{sect:space} that
$\RatCurves(X)$ is constructed as the quotient of the {\em
normalization} of $\Hom_{\operatorname{bir}} (\P^1, X)$. While
$\Hom_{\operatorname{bir}} (\P^1, X)$ is in 1:1-correspondence with
morphisms, $\P^1\to X$, that are birational onto their imnage, the
normalization morphism
$$
\Hombir (\P^1, X) \to \Hom_{\operatorname{bir}} (\P^1, X),
$$
need not be injective. For complex manifolds, however, we have the
following workaround.

\begin{lem}\label{lem:onetoone}
  Under the assumptions of Theorem~\ref{thm:irred}, let $x\in X$ be a
  general point and set
  $$
  H_x := \{\ell \in H \,|\, x \in \ell \}
  $$
  Then the closed points of $H_x$ are in 1:1-correspondence with
  the associated curves in $X$.
\end{lem}
\begin{proof}
  Since $x$ is a general point and since we have picked a fixed
  component, $H'$, all rational curves through $x$ are free by the
  proof of \cite[thm.~1.1]{KMM92}. The space
  $\Hom_{\operatorname{bir}} (\P^1, X)$ is therefore smooth at every
  point $f \in \Hom_{\operatorname{bir}} (\P^1, X)$ whose image
  contains the point $x$ by \cite[thm.~II.1.7]{K96}. The normalization
  morphism
  $$
  \Hombir (\P^1, X) \to \Hom_{\operatorname{bir}} (\P^1, X),
  $$
  is thus isomorphic in a neighborhood of $f$. Since
  $\Hom_{\operatorname{bir}} (\P^1, X)$ is in 1:1 correspondence with
  morphisms $\P^1\to X$, the claim follows.
\end{proof}

This enables us to prove Theorem~\ref{thm:irred}.

\begin{proof}[Proof of Theorem~\ref{thm:irred}]
  Choose a general point $x\in X$, and let $\tau : H \dasharrow
  \P(T_X^\vee)$ be the tangent morphism described in the
  introduction. Since all curves associated with $H_x$ are smooth,
  $\tau$ restricts to a regular morphism
  $$
  \tau_x : H_x \to \P(T_X^\vee|_x).
  $$ This morphism is known to be finite \cite[thm.~3.4]{Keb00a}. By
  Theorem~\ref{thm:main-criterion1}, $\tau_x$ is injective.
  
  Now assume that $H_x$ is not irreducible, $H_x = H_{x,1} \cup \ldots
  \cup H_{x,n}$. Since $\tau_x$ is finite, we have that
  $$
  \dim (\tau_x(H_{x,1}))+\dim (\tau_x(H_{x,2})) \geq \dim X-1 =
  \dim \P(T_X^\vee|_x)
  $$
  Thus, by \cite[thm.~I.7.2]{Ha77}, 
  $$
  \tau_x(H_{x,1}) \cap \tau_x(H_{x,2}) \not = \emptyset.
  $$
  It follows that $\tau$ is not injective, a contradiction.
\end{proof}

Lemma~\ref{lem:onetoone} raises the following interesting question.

\begin{question}
  Are there other conditions than smoothness over $\mathbb C$ which
  guarantee that closed points of $H_x$ are in 1:1-correspondence with
  rational curves?
\end{question}

\subsection{Stability of the tangent bundle, Automorphism Groups}

The setup of Theorem~\ref{thm:irred} naturally generalizes the notion
of a prime Fano manifold, i.e., one that is covered by lines under a
suitable embedding. Virtually all results that have been obtained for
prime Fanos hold in the more general setup of
Theorem~\ref{thm:irred}. We give two examples.

\begin{thm}
  In the setup of Theorem~\ref{thm:irred}, the tangent bundle $T_X$ is
  stable. In particular, $T_X$ is simple, i.e., the only endomorphisms
  of $T_X$ are the scalar multiplications.
\end{thm}
\begin{proof}
  The proof of \cite[thm.~2.11]{Hwa} applies verbatim.
\end{proof}

For any complex variety $X$, let $\Aut_0(X)$ denote the maximal
connected subgroup of the group of automorphisms. By universal
properties, an automorphism of a complex variety induces an
automorphism of the space $\RatCurves(X)$. It might be interesting to
note that in our setup the converse also holds.

\begin{thm}
  In the setup of Theorem~\ref{thm:irred}, the groups $\Aut_0(X)$ and
  $\Aut_0(H)$ coincide.
\end{thm}
\begin{proof}
  The proof of \cite[cor.~1]{HM02} applies verbatim.
\end{proof}

\subsection{Contact Manifolds}

Let $X$ be a projective contact manifold over $\mathbb C$, e.g.~the
twistor space over a Riemannian manifold with Quaternionic-K\"ahlerian
holonomy group and positive curvature. We refer to \cite{Kebekus01}
and the references therein for an introduction and for the relevant
background information.

If $X$ is different from the projective space, it has been shown in
\cite{Kebekus01} that $X$ is covered by a compact family of rational
curves $H \subset \RatCurves(X)$ such that for a general point $x$,
all curves associated with points in $H_x$ are smooth. Thus, the
assumptions of Theorem~\ref{thm:main-criterion2} are satisfied, and
$\tau$ is generically injective. This has been shown previously in
\cite{Keb00b} using rather involved arguments which heavily rely on
obstructions to deformations coming from contact geometry.

\bibliographystyle{alpha}

\end{document}